\pdfoutput=1
\documentclass[11pt]{article}
\usepackage[T1]{fontenc}
\usepackage{lmodern,microtype,amsmath,amssymb,amsthm,mathtools,tikz}
\usepackage[a4paper,margin=26mm]{geometry}
\usepackage{xurl}
\usepackage[colorlinks=true,linkcolor=blue!55!black,citecolor=blue!55!black,urlcolor=blue!55!black]{hyperref}
\hypersetup{pdftitle={Extinction times for synchronous majority dynamics on the square lattice},
  pdfauthor={Zifan Yang and Shuqing Zhong},
  pdfkeywords={Probabilistic cellular automata, interacting particle systems, extinction time, zero-temperature Glauber dynamics, censoring inequality}}

\newtheorem{theorem}{Theorem}[section]

\newtheorem{lemma}[theorem]{Lemma}
\newtheorem{corollary}[theorem]{Corollary}
\theoremstyle{remark}
\newtheorem{remark}[theorem]{Remark}

\newcommand{\PP}{\mathbb{P}}
\newcommand{\EE}{\mathbb{E}}
\newcommand{\Z}{\mathbb{Z}}
\newcommand{\one}{\mathbf{1}}
\newcommand{\floor}[1]{\lfloor #1\rfloor}
\newcommand{\ceil}[1]{\lceil #1\rceil}

\title{Extinction times for synchronous majority dynamics on the square lattice}
\author{
  Zifan Yang\thanks{These authors contributed equally to this work.}\textsuperscript{,}\thanks{School of Future Technologies, Beijing Institute of Technology, Beijing,
100081, China. Email: 
\texttt{yangzf@bit.edu.cn}, \texttt{sqzhong@bit.edu.cn}.}
  \and
  Shuqing Zhong\footnotemark[1]\textsuperscript{,}\footnotemark[2]\textsuperscript{,}\thanks{Corresponding author.}
}
\date{}

\begin{document}
\maketitle

\begin{abstract}
For synchronous four-neighbour majority dynamics with fair tie-breaking on \(\Z^2\), we prove a two-sided comparison between the extinction-time tails under parallel updates and under independent rate-one Poisson clocks. The comparison uses alternating-colour scans, the Peres--Winkler censoring inequality, and two embeddings of the update schedules. It transfers the continuous-time extinction window of Caputo, Martinelli, Simenhaus, and Toninelli to the synchronous process: for constants \(0<A<B\) and \(K,d>0\), the extinction time \(T_N\) of an initially occupied \(N\times N\) square belongs to \((AN^2,BN^2]\) with probability at least \(1-Ke^{-dN}\) for all sufficiently large \(N\). The filled square maximises the recovery time in stochastic order among finite perturbations fitting in a square of side \(N\). Its mean is of order \(N^2\), confirming a conjecture of Fat\`es, Marcovici, and Taati, and \(\EE T_N^r=\Theta(N^{2r})\) for every fixed \(r>0\). The exponential moments of \(T_N/N^2\) are also uniformly bounded as \(N\to\infty\) for every fixed positive parameter. We bound the probability of extinction at subquadratic times. For the filled square, the earliest possible extinction round is \(N\) for odd \(N\) and \(N-1\) for even \(N\), and the probability of extinction in that round is \(\exp[-\Theta(N^2)]\) as \(N\to\infty\).
\end{abstract}



\section{Introduction and main results}\label{sec:intro}

Replacing single-site continuous-time updates by simultaneous parallel updates can substantially alter the microscopic evolution of an interacting particle system. Whereas continuous-time Glauber dynamics updates one spin at a time and naturally admits reversible Markov generator techniques, a synchronous probabilistic cellular automaton updates all sites at once: in each round, adjacent spins consult each other's \emph{previous} state, creating local synchronisation effects that obstruct direct geometric comparisons.

Let \(Q_N = \{0,\ldots,N-1\}^2\). At each integer time \(k \ge 0\), every site \(v \in \Z^2\) carries a state in \(\{0,1\}\). Preassign independent fair coins to all site--round pairs. If \(j\) of the four nearest neighbours of \(v\) are in state \(1\) at time \(k\), its state in round \(k+1\) is \(0\) if \(j<2\), \(1\) if \(j>2\), and the coin for \((v,k+1)\) if \(j=2\). The site itself has no vote. For a deterministic configuration \(x\) that agrees with a background state \(b\in\{0,1\}\) outside a finite set, let \(T_x\) be the first round at which every site has state \(b\), with \(T_x=\infty\) if there is no such round. Define \(n(x)\) as the least \(N\) such that a translate of \(Q_N\) contains \(\{v:x(v)\ne b\}\), with \(n(x)=0\) when \(x\) is homogeneous. We investigate the worst-case expected recovery time and the filled-square extinction time:
\[
  \tau(N) = \sup_{n(x)=N} \EE T_x, \qquad T_N = T_{\one_{Q_N}}.
\]
The supremum is over deterministic initial configurations of exact size \(N\) in either background, with expectations taken in \([0,\infty]\). Here \(n(x)\) counts lattice sites along one side of the enclosing square. This is the diameter convention of~\cite[Section~1.2]{fmt}; even a two-site perturbation can have arbitrarily large \(n(x)\).

Fat\`es, Marcovici, and Taati called this rule \emph{Maj-Random-If-Equal}. They proved \(\tau(N)=O(N^3)\) and conjectured quadratic scaling~\cite[Proposition~4.1 and Question~7.4]{fmt}. For the corresponding continuous-time zero-temperature Ising dynamics, Fontes, Schonmann, and Sidoravicius~\cite[Theorem~1.3]{fss} proved an exponentially small probability of survival beyond a sufficiently large multiple of \(N^2\). Caputo, Martinelli, Simenhaus, and Toninelli supplied the matching lower tail and stated the two-sided window~\cite[Theorem~2 and Remark~2]{cmst}. The two update schedules behave differently even on a filled \(2\times2\) square. Each of the two top sites is occupied after a parallel round with probability \(1/2\). If instead only these two sites are updated, left then right, the right site remains occupied only when both coins equal \(1\), with probability \(1/4\).

Our first result confirms the quadratic conjecture of~\cite{fmt}. It gives a fixed quadratic-scale interval for the filled-square extinction time, with exponentially small probability outside that interval, and identifies the filled square as the worst initial configuration.

\begin{theorem}[Extinction on the quadratic scale]\label{thm:main}
There exist constants \(0 < A < B < \infty\) and \(K, d > 0\) such that for all sufficiently large \(N\),
\begin{equation}\label{eq:window}
  \PP(T_N \le AN^2) + \PP(T_N > BN^2) \le K e^{-dN}.
\end{equation}
Moreover, for every \(N \ge 1\) and every configuration \(x\) with \(n(x)=N\), the recovery time \(T_x\) is stochastically dominated by \(T_N\). In particular, \(\tau(N) = \EE T_N = \Theta(N^2)\) as \(N \to \infty\), and \(\EE T_N^r = \Theta(N^{2r})\) for every fixed \(r > 0\). For every fixed \(\lambda>0\), we also have
\[
  \limsup_{N\to\infty}\ \sup_{n(x)=N}\EE\bigl[e^{\lambda T_x/N^2}\bigr]<\infty.
\]
\end{theorem}

The lower-tail bound in~\eqref{eq:window} concerns the filled square. Two exceptional sites at lattice distance at least three, for example, vanish in one round. For the filled square, the next theorem identifies the earliest possible extinction round and the order of the logarithm of its probability.

\begin{theorem}[First possible extinction]\label{thm:atom}
For every \(N \ge 1\), set \(m_N = 1 + 2\floor{(N-1)/2}\) (so \(m_N = N\) for odd \(N\) and \(m_N = N-1\) for even \(N\)). Then \(T_N \ge m_N\) almost surely, and the atom \(q_N = \PP(T_N = m_N)\) satisfies
\[
  q_N \ge 2^{-2N^2}, \qquad
  q_{2h} \le 2^{-4h^2}, \qquad
  q_{2h+1} \le 4 \cdot 2^{-h^2} \qquad (h \ge 1).
\]
In particular, \(m_N = \min\{k \ge 0 : \PP(T_N = k) > 0\}\) and \(-\log q_N = \Theta(N^2)\) as \(N \to \infty\).
\end{theorem}

Fat\`es, Marcovici, and Taati already used a space--time parity split in their analysis of a boundary row~\cite[proof of Proposition~4.1]{fmt}. Peres and Winkler compared alternating and random update schemes through censoring~\cite[Section~3]{pw}, and dependency-depth bounds for Poisson update schedules are available in~\cite[Lemma~3.11]{ly}. The comparison here uses local readiness to extract complete alternating scans from a clock schedule, and uses the ranks of clock events to embed the schedule in the opposite direction. This gives two survival inequalities for the original synchronous process, with times differing by constant factors and explicit schedule error terms.

\paragraph{Overview of the proof.}
Section~\ref{sec:comp} proves the comparison and Section~\ref{sec:tails} applies it to the extinction window, moments and early-extinction bounds. The proof of Theorem~\ref{thm:atom} in Section~\ref{sec:early} uses only the discrete-time rule: the distance to the nearest corner determines when a site can first become zero and which coins must be zero for this to happen.

\section{A two-sided comparison of update schedules}\label{sec:comp}

By spin-flip symmetry (complementing all tie-breaking coins preserves the law) and translation invariance, we may assume without loss of generality that the background state is \(0\) and that the initial support \(\{v : x(v) = 1\}\) is contained in \(Q_N\). Every site outside \(Q_N\) has at most one nearest neighbour in \(Q_N\), so induction over rounds shows that all exterior sites remain \(0\) at every time \(k \ge 0\). Furthermore, for any fixed realisation of the tie-breaking coins, the local majority rule is monotone increasing in the input configuration. Coupling two trajectories via the same site-and-round coins therefore preserves inclusion of the occupied sets, which yields the stochastic domination \(T_x \preceq T_N\) in Theorem~\ref{thm:main} (and hence \(\tau(N) = \EE T_N\)).

To relate the synchronous dynamics to continuous time, consider the same zero-temperature majority rule on \(Q_N\) with all-zero exterior, driven by independent rate-one Poisson clocks at the sites of \(Q_N\) and a fresh independent fair coin at each clock ring. Starting both processes from the filled square \(\one_{Q_N}\), define the survival probabilities
\[
  p_N(k) = \PP(T_N > k), \qquad C_N(t) = \PP(T_N^{\mathrm{cl}} > t),
\]
where \(T_N^{\mathrm{cl}}\) is the first physical time at which the clock process reaches the all-zero configuration. Under either update schedule, the all-zero configuration is absorbing. Almost surely, on every bounded time interval there are only finitely many rings in \(Q_N\), and no two clocks ring simultaneously.

\begin{theorem}[Comparison of survival probabilities]\label{thm:bridge}
For all integers \(N \ge 2\) and \(k \ge 0\), set \(t_k = 8(\log 2)(k+1)\) and \(E_N(k) = N^2 2^{-(k+1)}\). Then
\begin{align}
  C_N(t_k) &\le p_N(k) + E_N(k), \label{eq:forward}\\
  p_N(256(k+1)+1) &\le 2 C_N(t_k) + 2 E_N(k). \label{eq:reverse}
\end{align}
\end{theorem}

The proof of Theorem~\ref{thm:bridge} rests on three lemmas established in the subsections below.

\subsection{Alternating-colour scans and censoring}

Colour each site \((i,j) \in Q_N\) \emph{even} or \emph{odd} according to the parity of \(i+j\). A \emph{scan} of a given colour updates all sites of that colour once; because sites of the same colour are pairwise nonadjacent, their updates commute. Let \(A_c(s)\) denote the probability that the occupied set is nonempty after \(s\) alternating scans starting with colour \(c \in \{\mathrm{even}, \mathrm{odd}\}\), initialised from \(\one_{Q_N}\).

For any finite sequence (word) \(w = (v_1,\ldots,v_m)\) of sites in \(Q_N\) and any function \(f : \{0,1\}^{Q_N} \to \mathbb{R}\), let \(K_w f\) denote the expectation of \(f\) at the configuration obtained by starting from \(\one_{Q_N}\) and applying the single-site majority updates along \(w\) in order, averaging over independent fair coins.

\begin{lemma}[Zero-temperature censoring and parity comparison]\label{lem:parity}
\leavevmode
\begin{enumerate}
\item[(a)] If a word \(u\) is obtained from \(w\) by deleting any subset of updates, then
\begin{equation}\label{eq:censor}
  K_w f \le K_u f \qquad \text{for every increasing function } f : \{0,1\}^{Q_N} \to \mathbb{R}.
\end{equation}
\item[(b)] For all integers \(k, q \ge 0\) and each starting colour \(c \in \{\mathrm{even},\mathrm{odd}\}\),
\begin{equation}\label{eq:parity}
  A_c(k+1) \le p_N(k), \qquad p_N(2q+1) \le 2 A_{\mathrm{even}}(2q).
\end{equation}
\end{enumerate}
\end{lemma}

\begin{proof}
(a) At finite inverse temperature \(\beta \in (0,\infty)\), inequality~\eqref{eq:censor} is the Peres--Winkler censoring theorem for the ferromagnetic Ising model started from its maximal configuration~\cite{pw}. In spin variables \(\sigma=2x-1\), the fixed minus neighbours contribute the one-site field \(h_v=-\#\{u\notin Q_N:u\sim_{\Z^2}v\}\), where \(\sim_{\Z^2}\) denotes nearest-neighbour adjacency on the lattice. The interior couplings remain ferromagnetic, so the finite-volume Gibbs weight satisfies the FKG lattice condition and its single-site conditionals are increasing in the other spins. The finite-volume Gibbs measure has full support at finite \(\beta\). The heat-bath probability of updating a site to \(+1\) given the total neighbour spin sum \(s\), including the fixed exterior spins, is \((1 + e^{-2\beta s})^{-1}\). As \(\beta \to \infty\), this converges to the zero-temperature majority kernel, including probability \(1/2\) at \(s=0\). Because the state space \(\{0,1\}^{Q_N}\) and both words \(u, w\) are finite, matrix multiplication of single-site transition kernels commutes with the limit \(\beta \to \infty\), yielding~\eqref{eq:censor}. In what follows, we apply~\eqref{eq:censor} conditionally on unmarked update schedules, which are independent of the tie-breaking coins.

(b) Run two alternating-scan processes from \(\one_{Q_N}\) with opposite initial colours. At slot \(s\), a site updated by a scan uses the coin of that site in synchronous round \(s\); the two scans use disjoint site--round coordinates of the synchronous coin field. By induction on the slot index \(s \ge 1\), immediately after slot \(s\) the active colour matches the synchronous configuration at round \(s\), while the inactive colour matches the synchronous configuration at round \(s-1\); indeed, every neighbour of an active site is inactive and therefore holds its round-\((s-1)\) value at the moment of the update.

Consequently, if the synchronous process has reached the all-zero state by round \(k\), then it is all-zero at both rounds \(k\) and \(k+1\), so either scan sequence is empty after \(k+1\) scans; this gives \(A_c(k+1) \le p_N(k)\). Conversely, if the synchronous process is still nonempty at round \(2q+1\), then at least one of the two scan sequences is nonempty after \(2q+1\) scans. Deleting the final scan from the even-start sequence and the first scan from the odd-start sequence leaves, in both cases, exactly \(2q\) alternating scans beginning with the even colour. Applying~\eqref{eq:censor} to the nonempty-set indicator and taking a union bound over the two sequences yields \(p_N(2q+1) \le 2 A_{\mathrm{even}}(2q)\).
\end{proof}

\subsection{Selecting complete scans from the Poisson clocks}

Write \(n_v(s)\) for the time of the first ring of the Poisson clock at site \(v \in Q_N\) strictly after time \(s\), and write \(u \sim v\) when \(u\) and \(v\) are nearest neighbours in \(Q_N\). On the common probability space of independent rate-one Poisson clocks, define recursively
\[
  D_0(v) = 0, \qquad
  R_j(v) = \max_{u \sim v} D_j(u), \qquad
  D_{j+1}(v) = n_v\bigl(R_j(v)\bigr) \qquad (j \ge 0,\ v \in Q_N).
\]

\begin{lemma}[Poisson-to-scan embedding]\label{lem:select}
For all integers \(N \ge 2\), \(k \ge 0\), real \(t > 0\), and \(c \in \{\mathrm{even},\mathrm{odd}\}\),
\begin{equation}\label{eq:dirone}
  C_N(t) \le A_c(k+1) + N^2 8^{k+1} e^{-t/2}.
\end{equation}
\end{lemma}

\begin{proof}
At each layer \(j \in \{0,\ldots,k\}\), retain the clock ring at time \(D_{j+1}(v)\) only for sites \(v\) of the colour active in slot \(j+1\) (starting with colour \(c\) at slot \(1\)). Use the usual right-continuous filtration of the joint clock process. The maximum of finitely many almost surely finite stopping times is again such a stopping time, and the first ring of a given clock strictly after such a time is also an almost surely finite stopping time. Induction therefore gives this property for \(R_j(v)\) and \(D_{j+1}(v)\). By the strong Markov property, the joint clocks after \(R_j(v)\) have fresh Poisson increments independent of the history through \(R_j(v)\); in particular, \(W_j(v)=D_{j+1}(v)-R_j(v)\) is \(\mathrm{Exp}(1)\)-distributed and independent of that history.

Let \(M = \max_{0 \le j \le k} \max_{v \text{ active in slot } j+1} D_{j+1}(v)\) denote the latest retained ring time across the first \(k+1\) layers. Since \(N \ge 2\), every site \(u \in Q_N\) has at least one neighbour \(v \in Q_N\); because \(u\) and \(v\) have opposite colours, if \(u\) is active in slot \(j+1\) then \(v\) is active in slot \(j+2\), and \(D_{j+1}(u) \le R_{j+1}(v) < D_{j+2}(v)\). Hence every retained ring in layers \(0,\ldots,k-1\) is strictly followed by a neighbour's retained ring in the next layer, so \(M\) is attained in the final layer \(j=k\). Expanding the recurrence along nearest-neighbour paths \((v_0, v_1, \ldots, v_k)\) in \(Q_N\) expresses \(M\) as the maximum of \(\sum_{j=0}^{k} W_j(v_j)\) over paths ending at a site of the active colour in the final layer. Along any \emph{fixed} path, \(D_{j+1}(v_j) \le R_{j+1}(v_{j+1})\), so \(W_0(v_0),\ldots,W_j(v_j)\) are measurable with respect to the history at \(R_{j+1}(v_{j+1})\). Conditioning successively from \(j=k\) down to \(j=0\) gives \(\EE[\exp(\frac{1}{2}\sum_{j=0}^{k} W_j(v_j))] = 2^{k+1}\). Since there are at most \(N^2 4^k\) nearest-neighbour paths of length \(k\) in \(Q_N\), Markov's inequality and a union bound yield the following estimate, without requiring independence across paths:
\begin{equation}\label{eq:selection}
  \PP(M > t) \le N^2 4^k 2^{k+1} e^{-t/2} \le N^2 8^{k+1} e^{-t/2}.
\end{equation}

On the event \(\{M \le t\}\), the retained rings all occur by time \(t\), with one update at every active site in each of the \(k+1\) layers. We show that they can be reordered into complete alternating scans without changing their endpoint. Call two rings \emph{conflicting} if their sites coincide or are nearest neighbours. By construction, neighbouring sites in consecutive layers satisfy \(D_j(u) < D_{j+1}(v)\); traversing an edge and returning shows \(D_j(v) < D_{j+2}(v)\). Consequently, distinct retained updates use distinct clock rings, and every conflicting pair occurs in chronological order consistent with its layer order, even across nonconsecutive layers. Any two linear orderings of a finite set of rings that agree on all conflicting pairs differ by a sequence of adjacent transpositions of nonconflicting rings (iteratively bubble the first desired ring to the front). Since nonconflicting updates commute when each retained ring carries its own coin, reordering the retained rings by layer preserves the endpoint configuration. Since the selection depends only on the unmarked clocks, the retained coins are independent fair coins conditional on the schedule. Deleting all unretained rings from the chronological clock word up to time \(t\) and applying Lemma~\ref{lem:parity}(a) conditionally on the unmarked clock schedule gives
\[
  C_N(t) \le A_c(k+1) + \PP(M > t),
\]
which combined with~\eqref{eq:selection} proves~\eqref{eq:dirone}.
\end{proof}

\subsection{Completing the Poisson schedule to alternating scans}

For the reverse direction, we use a dependency-depth estimate for the Poisson schedule (cf.~\cite[Lemma~3.11]{ly}) and pad each layer to a full scan.

\begin{lemma}[Scan-to-Poisson embedding]\label{lem:complete}
For all integers \(N \ge 2\), \(q \ge 1\), and real \(t > 0\),
\begin{equation}\label{eq:dirtwo}
  A_{\mathrm{even}}(2q) \le C_N(t) + \frac{1}{5} N^2 e^{9t} 2^{-q}.
\end{equation}
\end{lemma}

\begin{proof}
Enumerate all clock rings in \(Q_N\) during \((0,t]\) in chronological order, and recursively assign to each ring a \emph{rank} equal to \(1\) plus the maximum rank of all strictly earlier conflicting rings (with the empty maximum set to \(0\)). If some ring attains rank at least \(q \ge 1\), tracing predecessors back yields \(q\) chronologically ordered rings whose sites \(v_1,\ldots,v_q\) satisfy either \(v_{i+1} = v_i\) or \(v_{i+1} \sim v_i\) for each \(i\).

For any fixed list \((v_1,\ldots,v_q)\) of this form, define \(H_0 = 0\) and \(H_i = n_{v_i}(H_{i-1})\) for \(1 \le i \le q\). By the strong Markov property, the increments \(H_i - H_{i-1}\) are independent \(\mathrm{Exp}(1)\) variables (even when the list revisits sites), so \(\EE[e^{-9 H_q}] = 10^{-q}\). Whenever the Poisson schedule realises the list in chronological order within \((0,t]\), its last ring occurs at or after \(H_q\), so \(H_q \le t\). Since there are at most \(N^2 5^{q-1}\) such lists (five choices at each step: the current site or one of its four neighbours), Markov's inequality and a union bound give
\begin{equation}\label{eq:rank}
  P_{\mathrm{bad}} := \PP(\text{some ring in }(0,t]\text{ has rank}\ge q)
  \le N^2 5^{q-1} e^{9t} 10^{-q} = \frac{1}{5} N^2 e^{9t} 2^{-q}.
\end{equation}

On the complementary event where all rings have rank strictly below \(q\), assign each even-site ring of rank \(r \in \{1,\ldots,q-1\}\) to slot \(2r-1\) and each odd-site ring of rank \(r\) to slot \(2r\). Because earlier conflicting rings have strictly smaller rank, sorting the rings by slot preserves the relative order of every conflicting pair and hence leaves the endpoint configuration unchanged for fixed coins. Moreover, no two rings in the same slot can share a site (they would conflict and thus have distinct ranks), so each slot contains distinct sites of a single colour. Inserting all missing sites into each of the \(2q\) slots (including any initially empty slots) produces \(2q\) complete alternating scans starting with the even colour. Keep each original coin attached to its ring and give each inserted update a fresh fair coin independent of the clocks and of all other coins. Conditional on the unmarked clock schedule, the padded word therefore has the alternating-scan endpoint law. Deleting the inserted updates recovers the sorted clock word, so Lemma~\ref{lem:parity}(a) implies \(A_{\mathrm{even}}(2q) \le C_N(t) + P_{\mathrm{bad}}\), and~\eqref{eq:rank} yields~\eqref{eq:dirtwo}.
\end{proof}

\begin{proof}[Proof of Theorem~\ref{thm:bridge}]
Choose \(t = t_k = 8(\log 2)(k+1)\) and \(q = 128(k+1)\). Substituting \(t_k\) into~\eqref{eq:dirone} gives an error term \(N^2 8^{k+1} 2^{-4(k+1)} = N^2 2^{-(k+1)} = E_N(k)\); combining this with the first inequality in~\eqref{eq:parity} yields~\eqref{eq:forward}. Similarly, substituting \(t_k\) and \(q\) into~\eqref{eq:dirtwo} gives an error term \(\tfrac{1}{5} N^2 2^{72(k+1)} 2^{-128(k+1)} = \tfrac{1}{5} N^2 2^{-56(k+1)} \le E_N(k)\); combining this with the second inequality in~\eqref{eq:parity} yields~\eqref{eq:reverse}.
\end{proof}

\begin{remark}[Scope of the comparison method]\label{rem:scope}
The schedule constructions extend to finite bipartite interaction graphs without isolated vertices, provided each single-site update depends on its neighbours and not on its own previous state. For a binary monotone heat-bath system started from its maximal configuration, the censoring inequality then gives the corresponding comparisons of scan and clock endpoint expectations, with path counts controlled by the maximum degree. The parity construction recovers the synchronous configuration from two scan trajectories. The survival comparison further requires that the all-zero configuration be absorbing. At positive temperature in the Ising model, the schedule embeddings and censoring comparisons remain available, but nonemptiness at a given time no longer identifies survival.
\end{remark}

\section{The extinction window, moments, and early-extinction tails}\label{sec:tails}

Caputo, Martinelli, Simenhaus, and Toninelli~\cite[Theorem~2]{cmst} analyse the continuous-time zero-temperature heat-bath Ising dynamics started from all minus spins on \(\{-L,\ldots,L\}^2\) with plus exterior, driven by one rate-one Poisson clock per site. Under the change of variables \(x = (1-\sigma)/2\) and translation by \((L,L)\), their process is identical in law to ours on \(Q_{2L+1}\), so their theorem provides constants \(0 < c \le 1\) and \(\gamma > 0\) such that
\begin{equation}\label{eq:cmst}
  \PP\bigl(c L^2 \le T_{2L+1}^{\mathrm{cl}} \le L^2/c\bigr) \ge 1 - e^{-\gamma L} \qquad (L \ge 1).
\end{equation}
In particular, the strict survival probability \(C_{2L+1}(t)\) is at least \(1 - e^{-\gamma L}\) for every \(t < c L^2\).

To accommodate even side lengths \(N\), observe that for \(u \le v\), restricting the clocks and coins of \(Q_v\) to \(Q_u\) (and extending the smaller filled-square process by \(0\) outside \(Q_u\)) preserves the product law of the restricted clocks and maintains inclusion of occupied sets at every ring by monotonicity. Induction over rings yields \(C_u(t) \le C_v(t)\) for all \(t \ge 0\). Setting \(L_- = \floor{(N-1)/2}\) and \(L_+ = \ceil{(N-1)/2}\), we obtain
\begin{equation}\label{eq:sandwich}
  C_{2L_- + 1}(t) \le C_N(t) \le C_{2L_+ + 1}(t) \qquad (t \ge 0).
\end{equation}

\begin{proof}[Proof of Theorem~\ref{thm:main}]
Write \(a = 8\log 2\). For all sufficiently large \(N\), define
\[
  k_- = \floor{\frac{c L_-^2}{2a}} - 1, \qquad
  k_+ = \ceil{\frac{L_+^2}{ca}} - 1, \qquad
  b_N = 256(k_+ + 1) + 1.
\]
For such \(N\) we have \(L_-\ge1\) and \(k_-,k_+\ge0\), as required in Theorem~\ref{thm:bridge} and in~\eqref{eq:cmst}. Then \(t_{k_-} = a(k_- + 1) \le \frac{1}{2} c L_-^2 < c L_-^2\) and \(t_{k_+} = a(k_+ + 1) \ge L_+^2/c\). Combining~\eqref{eq:forward}, \eqref{eq:reverse}, \eqref{eq:cmst}, and~\eqref{eq:sandwich} yields
\begin{align}
  1 - p_N(k_-) &\le e^{-\gamma L_-} + E_N(k_-), \label{eq:window1}\\
  p_N(b_N) &\le 2 e^{-\gamma L_+} + 2 E_N(k_+). \label{eq:window2}
\end{align}
Since \(L_\pm=N/2+O(1)\), both \(k_-+1\) and \(b_N\) are of order \(N^2\). Choose \(0<A<B<\infty\) so that \(AN^2\le k_-\) and \(b_N\le BN^2\) for all large \(N\). The comparison errors \(E_N(k_\pm)\) are \(O(N^2e^{-c'N^2})\) for some \(c'>0\), and \(e^{-\gamma L_\pm}=O(e^{-\gamma N/2})\). Thus~\eqref{eq:window1}--\eqref{eq:window2} imply~\eqref{eq:window}.

To control the upper tail beyond \(b_N\) and deduce the moment bounds, let \(\rho_N = p_N(b_N)\). At every deterministic block boundary \(j b_N\) (\(j \ge 0\)), the configuration remains supported in \(Q_N\) and the subsequent \(b_N\) rounds use fresh independent coins. The filled square dominates every state supported in \(Q_N\), including one whose remaining support has smaller diameter. Conditional on survival up to round \(j b_N\), stochastic domination by the filled square bounds the probability of surviving the next block by \(\rho_N\). Iterating over blocks gives, for every \(x\) with \(n(x) = N\) and every \(j \ge 0\),
\begin{equation}\label{eq:restart}
  \PP(T_x > j b_N) \le \rho_N^j.
\end{equation}
By~\eqref{eq:window2}, \(\rho_N=O(e^{-d' N})\) for some \(d'>0\), and in particular \(\rho_N\le1/2\) for all sufficiently large \(N\). The tail-sum identity for nonnegative integer-valued extended recovery times, grouped into blocks of length \(b_N\), gives
\[
  \EE T_x=\sum_{k\ge0}\PP(T_x>k)
  \le b_N\sum_{j\ge0}\rho_N^j
  =\frac{b_N}{1-\rho_N}\le 2b_N.
\]
The lower bound \(\EE T_N \ge k_- p_N(k_-)\) follows from~\eqref{eq:window1}; hence \(\tau(N) = \EE T_N = \Theta(N^2)\).

The block estimate also controls exponential moments on the quadratic scale. Choose \(B_0>0\) with \(b_N\le B_0 N^2\) for all sufficiently large \(N\), and fix \(\lambda>0\). Because \(\rho_N\to0\), we have \(\rho_N e^{\lambda B_0}\le1/2\) for all sufficiently large \(N\), depending on \(\lambda\). The block tail implies \(T_x<\infty\) almost surely. Decomposing its finite values into \(T_x=0\) and \(j b_N<T_x\le(j+1)b_N\), and using~\eqref{eq:restart}, we get
\[
  \sup_{n(x)=N}\EE\bigl[e^{\lambda T_x/N^2}\bigr]
  \le 1+\sum_{j\ge0}e^{\lambda B_0(j+1)}\rho_N^j
  =1+\frac{e^{\lambda B_0}}{1-\rho_N e^{\lambda B_0}}
  \le 1+2e^{\lambda B_0}.
\]
This proves the asserted uniform bound for every fixed positive parameter. As \(z^r\le C_{r,\lambda} e^{\lambda z}\) for all \(z\ge0\), it gives \(\EE T_N^r=O(N^{2r})\). The lower tail~\eqref{eq:window1} gives the matching bound \(\EE T_N^r\ge k_-^r p_N(k_-)=\Omega(N^{2r})\).
\end{proof}

Packing disjoint sub-squares into \(Q_N\) amplifies the lower-tail estimate of Theorem~\ref{thm:main} to arbitrary sub-quadratic time scales.

\begin{corollary}[Early-extinction bounds]\label{cor:packing}
There exists a constant \(c_1 > 0\) such that for every integer \(N \ge 1\) and every integer \(k \ge 0\),
\[
  \PP(T_N \le k) \le 2 \exp\!\left(-\frac{c_1 N^2}{k+1}\right).
\]
Moreover, there are constants \(a_0,c_2>0\) such that whenever \(k+1\le a_0 N^2\),
\[
  \PP(T_N \le k) \le \exp\!\left(-\frac{c_2 N^2}{\sqrt{k+1}}\right).
\]
\end{corollary}

\begin{proof}
Write \(F_N(k)=1-p_N(k)\). Choose an integer \(D\ge1\) so that~\eqref{eq:window} holds for every side length at least \(D\), \(AD^2\ge1\), \(K e^{-dD}\le1/2\), and \(K\le e^{dD/2}\). For \(m=\ceil{D\sqrt{k+1}}\), we have \(m\ge D\) and \(k<Am^2\), hence \(F_m(k)\le\min\{1/2,e^{-dm/2}\}\).

Set \(s=\floor{N/m}\) and place \(s^2\) disjoint translates of \(Q_m\) in \(Q_N\). Evolve each sub-square with its own exterior fixed at \(0\), using the coins of the larger process at the sites it updates. The sub-square processes read disjoint sets of coins, so they are independent. Each is dominated by the restriction of the larger process to that sub-square: clamping its exterior to \(0\) can only decrease the outputs of its monotone updates. Thus extinction of the larger process by round \(k\) requires extinction of every sub-square by round \(k\), and \(F_N(k)\le F_m(k)^{s^2}\) if \(s\ge1\).

When \(N\ge m\), we have \(s\ge N/(2m)\) and \(m\le(D+1)\sqrt{k+1}\). The bound \(F_m(k)\le1/2\) gives the first estimate with \(c_1=(\log 2)/(4(D+1)^2)\). For \(N<m\) its right-hand side is greater than \(1\), so it holds as well. If \(k+1\le N^2/(D+1)^2\), then \(m\le N\), and the bound \(F_m(k)\le e^{-dm/2}\) gives
\[
  F_N(k)\le e^{-dm s^2/2}\le
  \exp\!\left(-\frac{dN^2}{8(D+1)\sqrt{k+1}}\right).
\]
This proves the second estimate with \(a_0=(D+1)^{-2}\) and \(c_2=d/(8(D+1))\).
\end{proof}

\section{The lower endpoint of the extinction law}\label{sec:early}

In this section we work entirely in discrete time to prove Theorem~\ref{thm:atom}. For each site \(v = (i,j) \in Q_N\), define its \emph{corner depth} (lattice distance to the nearest corner of \(Q_N\)) by
\[
  d(v) = \min(i, N-1-i) + \min(j, N-1-j).
\]
Note that the maximum depth in \(Q_N\) is \(2\floor{(N-1)/2} = m_N - 1\).

\begin{figure}[!htbp]
\centering
\begin{tikzpicture}[font=\small]
\foreach \side/\panelshift in {3/0,4/4.6}{
  \begin{scope}[xshift=\panelshift cm, x=.62cm, y=.62cm]
  \pgfmathtruncatemacro{\lastsite}{\side-1}
  \foreach \row in {0,...,\lastsite}{
    \foreach \col in {0,...,\lastsite}{
      \pgfmathtruncatemacro{\firstzero}{1+min(\col,\lastsite-\col)+min(\row,\lastsite-\row)}
      \pgfmathtruncatemacro{\tone}{8+12*(\firstzero-1)}
      \filldraw[fill=black!\tone,draw=black!65] (\col,\row) rectangle ++(1,1);
      \node at (\col+.5,\row+.5) {\firstzero};
    }
  }
  \node[below] at (.5*\side,-.18) {\(N=\side\)};
  \end{scope}
}
\end{tikzpicture}
\caption{First zero round at each site when every tie resolves to \(0\). Each entry equals \(1+d(v)\), with maximum \(3\) in both squares (attained at one site in each \(2\times2\) corner square for \(N=4\)). TikZ code assisted by OpenAI Codex (desktop build 26.1117.4957.0) and verified by the authors.}\label{fig:depth}
\end{figure}

\begin{proof}[Proof of Theorem~\ref{thm:atom}]
We first establish the deterministic lower bound \(T_N \ge m_N\) and the lower estimate on \(q_N\).

\smallskip
\noindent\emph{Step 1: Deterministic protection and the all-zero tie trajectory.}
We claim that for every realisation of the tie-breaking coins and every round \(t \ge 0\), all sites \(v \in Q_N\) with \(d(v) \ge t\) are occupied at round \(t\). For \(t = 0\) this is the initial condition, and for \(N = 1\) the induction step is vacuous. Assuming the claim at round \(t\), take \(v \in Q_N\) with \(d(v) \ge t+1\). Then \(v\) is not a corner of \(Q_N\), so at least three of its four nearest neighbours lie in \(Q_N\); since \(d\) changes by at most \(1\) along each edge, all such neighbours have depth at least \(t\) and are occupied at round \(t\). A strict majority therefore forces \(v\) to be occupied at round \(t+1\), completing the induction. Since the maximum depth in \(Q_N\) is \(m_N - 1\), at least one site remains occupied at round \(m_N - 1\), so \(T_N \ge m_N\) almost surely.

Next, consider the trajectory in which every tie resolves to \(0\). Its occupied set at round \(t\) is precisely \(\{v\in Q_N:d(v)\ge t\}\). For the reverse inclusion, suppose \(v\) is occupied at round \(t+1\) and that the assertion holds at \(t\). At least three neighbours of \(v\) were occupied at round \(t\), so they include an opposite pair along a coordinate axis. Both sites of this pair are in \(Q_N\). The neighbour toward a closest endpoint of that coordinate interval has depth \(d(v)-1\); if the two endpoints are equally close, either neighbour has that depth. The induction hypothesis gives \(d(v)-1\ge t\), and hence \(d(v)\ge t+1\). Thus the all-zero tie trajectory reaches extinction at round \(m_N=\max_v d(v)+1\) (Figure~\ref{fig:depth}).

Although the all-zero tie prescription specifies infinitely many coins, only two coins per site are actually needed to reproduce this trajectory. Specifically, set the coins at each site \(v \in Q_N\) at rounds \(d(v)+1\) and \(d(v)+2\) to \(0\). By Step~1, \(v\) is deterministically occupied prior to round \(d(v)+1\); at rounds \(d(v)+1\) and \(d(v)+2\), induction ensures the same neighbour inputs and tie outcome as on the all-zero schedule; and at any round \(s \ge d(v)+3\), every interior neighbour \(u \sim v\) has depth \(d(u) \le d(v)+1 < s-1\), so all neighbours of \(v\) were already \(0\) at round \(s-1\) and \(v\) remains \(0\) deterministically. Hence the two trajectories coincide, which gives \(q_N \ge 2^{-2N^2}\).

\smallskip
\noindent\emph{Step 2: Coins forced by extinction at round \(m_N\).}
Put \(h = \floor{N/2}\). In the southwest \(h \times h\) corner square \(\{0,\ldots,h-1\}^2\), every site \(v = (i,j)\) has depth \(d(v) = i + j\). Say that \(v\) \emph{vanishes on time} if it is in state \(0\) at round \(i+j+1\). At round \(i+j\), both its north neighbour \((i,j+1)\) and east neighbour \((i+1,j)\) lie in \(Q_N\) and have depth at least \(i+j\), so both are deterministically occupied by Step~1. For \(v\) to vanish on time, its west and south neighbours must therefore be \(0\) at round \(i+j\) \emph{and} its own tie-breaking coin at round \(i+j+1\) must be \(0\). In particular, any west or south neighbour of \(v\) that lies inside \(Q_N\) must also have vanished on time. Recursing toward the southwest corner \((0,0)\) forces one specified zero coin (at round \(i'+j'+1\)) at each of the \((i+1)(j+1)\) sites \((i',j') \in \{0,\ldots,i\} \times \{0,\ldots,j\}\). By lattice reflection symmetry, the analogous implication holds in each of the four \(h \times h\) corner squares.

For even \(N = 2h\) (\(h \ge 1\)), the four \(h \times h\) corner squares partition \(Q_N\), and their four tips farthest from their respective corners have depth \(2h-2\). If extinction occurs at \(m_{2h} = 2h-1\), all four tips must vanish on time. Since the four corner squares are pairwise disjoint, this forces one specified zero coin at each of the \(4h^2 = N^2\) distinct sites of \(Q_N\), yielding \(q_{2h} \le 2^{-4h^2}\).

For odd \(N = 2h+1\) (\(h \ge 1\)), the four tips of the corner squares are the four diagonal neighbours of the central site \((h,h)\), each of depth \(2h-2\). Suppose that all four tips are still occupied at round \(2h-1\). Since the centre has depth \(2h\), it is also occupied at round \(2h-1\). Each of the four axial neighbours of the centre then sees the centre and two occupied diagonal neighbours at round \(2h-1\), so all four axial neighbours are occupied at round \(2h\), which in turn forces the centre to survive round \(2h+1 = m_{2h+1}\). Consequently, extinction at round \(2h+1\) requires at least one of the four tips to vanish on time. Each such event forces \(h^2\) independent zero coins in its corresponding corner square, so a union bound gives \(q_{2h+1} \le 4 \cdot 2^{-h^2}\). This completes the proof of Theorem~\ref{thm:atom}.
\end{proof}

\begin{remark}
Combining the two parity bounds in Theorem~\ref{thm:atom} gives the uniform upper bound \(q_N \le 2^{-\floor{N^2/8}}\) for all \(N \ge 7\), since \(8(h^2 - 2) \ge (2h+1)^2\) for \(h \ge 3\).
\end{remark}

\section{Concluding remarks and open problems}\label{sec:remarks}

We conclude with several remarks and natural questions left open by our analysis.

\begin{enumerate}
\item \emph{Scaling limit of the extinction time.} For the continuous-time square droplet, \(T_N^{\mathrm{cl}}/N^2\) converges in probability to \(1/2\)~\cite[after Theorem~2.2]{lst}. For smooth droplets, Cerf and Louhichi~\cite{cl} computed the initial area drift near a boundary point in terms of its curvature and orientation. Lacoin, Simenhaus and Toninelli~\cite{lstgeneral} proved convergence of the evolving boundary to anisotropic curve-shortening flow for a class of smooth domains that includes nonconvex droplets. The two inequalities of Theorem~\ref{thm:bridge} compare different multiples of time, so they do not identify a limiting constant for \(T_N/N^2\). It remains open whether \(T_N/N^2\) converges in probability to a deterministic constant, or whether \(N^{-2}\EE T_N\) converges.
\item \emph{Exact asymptotics at the lower endpoint.} Theorem~\ref{thm:atom} determines the order \(-\log q_N = \Theta(N^2)\), leaving a constant-factor gap between the lower and upper exponents even along squares of a fixed parity. Identifying \(\lim_{h\to\infty} (2h)^{-2} \log q_{2h}\) and \(\lim_{h\to\infty} (2h+1)^{-2} \log q_{2h+1}\) remains open.
\item \emph{Non-bipartite lattices.} Our comparison in Section~\ref{sec:comp} makes essential use of the two-colouring of \(\Z^2\) to recover the two colour classes of a synchronous round from two alternating-scan processes with opposite starting colours. Extending the quadratic scaling analysis to synchronous majority rules on non-bipartite graphs---such as the triangular lattice---would require a different comparison mechanism.
\item \emph{Explicit constants.} The independent finite-volume argument in the companion Lean formalisation~\cite{artifact} gives
\[
  \frac{N^2}{850000}\le\tau(N)\le664066N^2\qquad(N\ge16).
\]
Its lower bound holds for every \(N\ge1\). The same development gives the asymptotic bounds
\[
  \frac{1}{1152000\log 2}\le
  \liminf_{N\to\infty}\frac{\EE T_N}{N^2}
  \le\limsup_{N\to\infty}\frac{\EE T_N}{N^2}\le331776.
\]
The last inequality follows from the archived estimate \(\tau(N)\le(331776 N^2+257)/(1-a_N)\), where \(0\le a_N\le5/N\) for \(N\ge16\). The proof of Theorem~\ref{thm:main} above transfers the sharper continuous-time probability window~\eqref{eq:cmst} and inherits its unspecified constants.
\end{enumerate}

\section*{Computational methods and formal verification}
The companion formalisation uses Lean~4.35.0-rc4~\cite{artifact}. OpenAI Codex (OpenAI, desktop build 26.1117.4957.0) assisted with exploring proof strategies and with writing and debugging the Lean code. The archived declaration \path{SquareRecoveryTime.two_way_clock_tail} uses the common error \(\tfrac65 N^2 2^{-(k+1)}\), whereas the coefficient \(1\) in Theorem~\ref{thm:bridge} follows here from the separate estimates in Lemmas~\ref{lem:select} and~\ref{lem:complete}. The literature-based reduction in \texttt{ClockLiteratureTransfer.lean} is formalised with the continuous-time window as a hypothesis; the archive also provides an independent finite-volume proof of quadratic recovery that does not take this window as a hypothesis. The full exponential window~\eqref{eq:window}, its strengthened moment consequence and the second estimate in Corollary~\ref{cor:packing} are established by the paper's arguments. Figure~\ref{fig:depth} evaluates the integer formula \(1+d(v)\) in TikZ, with code generation assisted by the same version of OpenAI Codex.

\section*{CRediT authorship contribution statement}
\textbf{Zifan Yang:} Conceptualization, Methodology, Formal analysis, Writing -- original draft.
\textbf{Shuqing Zhong:} Conceptualization, Methodology, Software (Lean formalisation), Writing -- review \& editing.

\section*{Funding}
This research did not receive any specific grant from funding agencies in the public, commercial, or not-for-profit sectors.

\section*{Declaration of competing interest}
The authors declare that they have no known competing financial interests or personal relationships that could have appeared to influence the work reported in this paper.

\section*{Data and code availability}
The Lean formalisation and accompanying finite-volume proof are publicly available at the fixed repository revision cited in~\cite{artifact} and are included in the supplementary archive \texttt{companion.zip}. For the comparison, see \texttt{SquareRecoveryTime.lean}; for the finite-volume quadratic bounds, see \texttt{AllScaleLower.lean} and \texttt{OriginalQuadraticUpper.lean}. The asymptotic estimates use \texttt{LowerConsequences.lean} and \texttt{SharpRecoveryTail.lean}. The first estimate in Corollary~\ref{cor:packing} is also formalised in \texttt{FastTail.lean}. The declarations in \texttt{EarliestSquare.lean}, \texttt{EarliestAtom.lean} and \texttt{MinimumAtomUpper.lean} establish the earliest extinction round and the two-sided bounds on its probability.

\section*{Declaration of generative AI and AI-assisted technologies in the manuscript preparation process}
During the preparation of this work, the authors used OpenAI Codex (OpenAI, desktop build 26.1117.4957.0) to improve exposition, readability, and LaTeX formatting, to obtain automated feedback on clarity, and to assist with the TikZ code for Figure~\ref{fig:depth}. Its use in exploring proof strategies and developing the Lean formalisation is described above. After using this tool, the authors reviewed and edited the content as needed and take full responsibility for the content of the published article.

\end{document}